\documentclass[journal]{IEEEtran}

\usepackage{xcolor,soul,framed} 
\usepackage{etoolbox}
\makeatletter
\patchcmd{\@makecaption}
  {\scshape}
  {}
  {}
  {}
  \makeatother
  
\colorlet{shadecolor}{yellow}
\usepackage[pdftex]{graphicx}
\graphicspath{{../pdf/}{../jpeg/}}
\DeclareGraphicsExtensions{.pdf,.jpeg,.png}

\usepackage[cmex10]{amsmath}
\usepackage{amssymb}
\ifCLASSOPTIONcompsoc
 \usepackage[caption=false,font=normalsize,labelfont=sf,textfont=sf]{subfig}
\else
 \usepackage[caption=false,font=footnotesize]{subfig}
\fi
\usepackage{algorithm,algorithmic}
\usepackage{array}
\usepackage{mdwmath}
\usepackage{mdwtab}
\usepackage{eqparbox}
\usepackage{url}

\usepackage[textwidth=1.5in]{todonotes}

\begin{document}
    \title{Koopman Operator Applications in Signalized Traffic Systems}
  \author{Esther~Ling,~\IEEEmembership{Student Member,~IEEE,}
      Liyuan~Zheng,~\IEEEmembership{Student Member,~IEEE,} \\
      Lillian~J.~Ratliff,~\IEEEmembership{Member,~IEEE,}
      and~Samuel~Coogan,~\IEEEmembership{Member,~IEEE}


  \thanks{This work was funded in part by the National Science Foundation under Awards 1749357 and 1736582}
  \thanks{E. Ling is with the School of Electrical and Computer Engineering, Georgia Institute of Technology, Atlanta, GA 30332 USA (e-mail: lingesther@gatech.edu).}
  \thanks{S. Coogan is with the School of Electrical and Computer Engineering and the School of Civil and Environmental Engineering, Georgia Institute of Technology, Atlanta, GA 30332 USA (e-mail: sam.coogan@gatech.edu).}
  \thanks{L. Zheng and L. Ratliff are with the Department of Electrical Engineering, University of Washington, Seattle, WA, 98195 USA (e-mail: liyuanz8@uw.edu; ratliffl@uw.edu).}}


\maketitle

\begin{abstract}

  This paper proposes Koopman operator theory and the related algorithm dynamical mode decomposition (DMD) for analysis and control of signalized traffic flow networks. DMD provides a model-free approach for representing complex oscillatory dynamics from measured data, and we study its application to several problems in signalized traffic. We first study a single signalized intersection, and we propose applying this method to infer traffic signal control parameters such as phase timing directly from traffic flow data. Next, we propose using the oscillatory modes of the Koopman operator, approximated with DMD, for early identification of unstable queue growth that has the potential to cause cascading congestion. Then we demonstrate how DMD can be coupled with knowledge of the traffic signal control status to determine traffic signal control parameters that are able to reduce queue lengths. Lastly, we demonstrate that DMD allows for determining the structure and the strength of interactions in a network of signalized intersections. All examples are demonstrated using a case study network instrumented with high resolution traffic flow sensors.
\end{abstract}

\begin{IEEEkeywords}
Koopman Operator, big data, time-series, dynamic mode decomposition, queue modeling, signal and phase estimation, instability detection, traffic prediction
\end{IEEEkeywords}

%
\IEEEpeerreviewmaketitle


\section{Introduction}

Traffic sensors at signalized intersections are becoming ubiquitous as more cities move towards enabling intelligent transportation systems. Developing real-time event monitoring and forecasting systems remains an open problem, with the large data streams presenting a rich source to mine. In this regard, traffic forecasting has been well studied, ranging from auto-regressive integrated moving average (ARIMA) models \cite{MIN2011606}, \cite{XMin_arima}, partial least-squares regression \cite{Coogan_PLS}, \cite{Dutreix_PLS}, neural networks \cite{Duan_NN}, \cite{Zheng_NN} and multi-variable state-space methods \cite{STATHOPOULOS_state}. Other than predicting future flow, a good forecasting model also provides insight into traffic patterns. Understanding the effect of propagation patterns can inform strategies on addressing bottlenecks and congestion. Some studies that focus on such qualitative analysis include \cite{Hofleitner2012}, \cite{Kumar_structureMatrix} and \cite{ding2018detecting}.
In this paper, we pursue both avenues by using a model that makes quantitative predictions as well as provides qualitative information on characteristics of the intersection.
A major difficulty in modeling traffic is its nonlinear behavior \cite{KAMARIANAKIS-regressions}. This nonlinearity is especially pronounced at signalized intersections due to large short-term fluctuations induced by traffic signaling and unstable conditions arising from congestion \cite{STATHOPOULOS_state}. According to \cite{Kamarianakis_regime}, higher oscillations in arterials decrease forecasting accuracy of models by 10\%-20\%. This suggests a need for developing new approaches that especially leverage high resolution measurements for signalized traffic.


Inspired by recent advances in data-driven modeling of nonlinear and oscillatory dynamical systems, we propose using the \emph{Koopman operator} framework for studying both quantitative and qualitiative properties of signalized traffic networks using high-resolution---in both space and time---data. The Koopman operator framework jointly captures the ideas of (i) quantitative prediction, by modeling the underlying nonlinear dynamics via a corresponding linear operator that exactly captures the dynamics but operates in an infinite dimensional space;
and (ii) qualitative understanding, by studying the spectral properties of (finite-dimensional approximations) of this linear operator. Postulated by Koopman in 1931 \cite{Koopman}, recent advances in applied Koopman Operator theory \cite{appliedKoopmanism}, \cite{kmd_power_2} have led to a renewed interest in the theory, with applications in the fields of power systems \cite{kmd_power_2}, \cite{susuki-power-coherency}, \cite{kmd_power_1}, \cite{raak-power-partition}, image processing \cite{brunton_dmdRPCA}, quantitative finance \cite{mann_trading}, disease modeling \cite{proctor-disease} and neuro-science \cite{BRUNTON_Neural}. This paper aims to leverage Koopman Operator theory and its  promising applications for large data in disparate domains to address quantitative and qualitative problems in signalized traffic.



The main contributions of this paper are as follows. First, we provide an analysis of the structure of the Koopman operator learned from traffic flow data using the numerical technique of \emph{dynamic mode decomposition}. We then show that, from this operator, we are able to extract Koopman modes that provide information about an intersection's timing plan and phase-splits. Next, we propose a Koopman-based instability analysis technique for automated detection of prolonged growth in traffic queues. Lastly, we highlight new analysis in modeling queues under modified timing parameters.


The remainder of the paper is organized as follows. Section \ref{sec:preliminaries} presents preliminary notation and definitions. In Section \ref{sec:dmd}, we provide brief theoretical background on the Koopman Operator and dynamic mode decomposition (DMD) and a review of Koopman Operator applications in other domains. In the next three sections, we highlight three applications of the Koopman Operator to signalized traffic at a single intersection: in Section \ref{sec:interpr-modes-sign}, we demonstrate Koopman Mode analysis for the problem of inferring signal and phase timing information from traffic flow data. In particular, we explain how to recover cycle times, phase sequence and green-splits from measured vehicle flows. In Section \ref{sec:instability-analysis}, we use Koopman Eigenvalues for identifying extended growth in queue dynamics. In Section~\ref{sec:learn-dynam-contr}, we model queue dynamics using signal phases as exogenous input using dynamic mode decomposition with control. We show the effect of varying green times on the reconstructed queue at a particular leg.  In Section \ref{sec:ident-oper-struct}, we consider a network of intersections and examine how the operator learned by DMD can be used to infer structural properties of the underlying network, drawing comparisons with ARIMA methods for multi-step prediction.  Finally, we provide concluding remarks and future directions in Section \ref{sec:conclusion}.

\section{Preliminaries}
\label{sec:preliminaries}
This paper considers networks of signalized intersections equipped with in-lane sensors that measure traffic flow. We will also sometimes assume direct access to queue measurements at the intersections, which are obtained by integrating the difference of the flow at upstream and downstream sensors.  We assume an intersection consists of up to four \emph{legs}: South-Bound (SB), North-Bound (NB), East-Bound (EB) and West-Bound (WB). Each leg has up to three \emph{movements}: Left-Turn (LT), Through (T) and Right-Turn (RT). Therefore, collectively, an intersection may have up to 12  turn-movements. A \emph{phase} consists of a collection of movements. Traffic flow at an intersection is governed by \emph{phase-splits} or \emph{green-splits} which allocate green time for each phase. The \emph{cycle time} is the total time allotted to complete all phases once. The \emph{phase sequence} is the relative ordering of each phase. Multiple phases can be active at once, and a \emph{barrier} separates a group of phases with non-conflicting movements that are allowed to be simultaneously active; for example, the phases for \{NBT, SBT\} might form a valid barrier.

\subsection{Dataset}
\begin{figure}
  \centering
\includegraphics[width=.9\columnwidth]{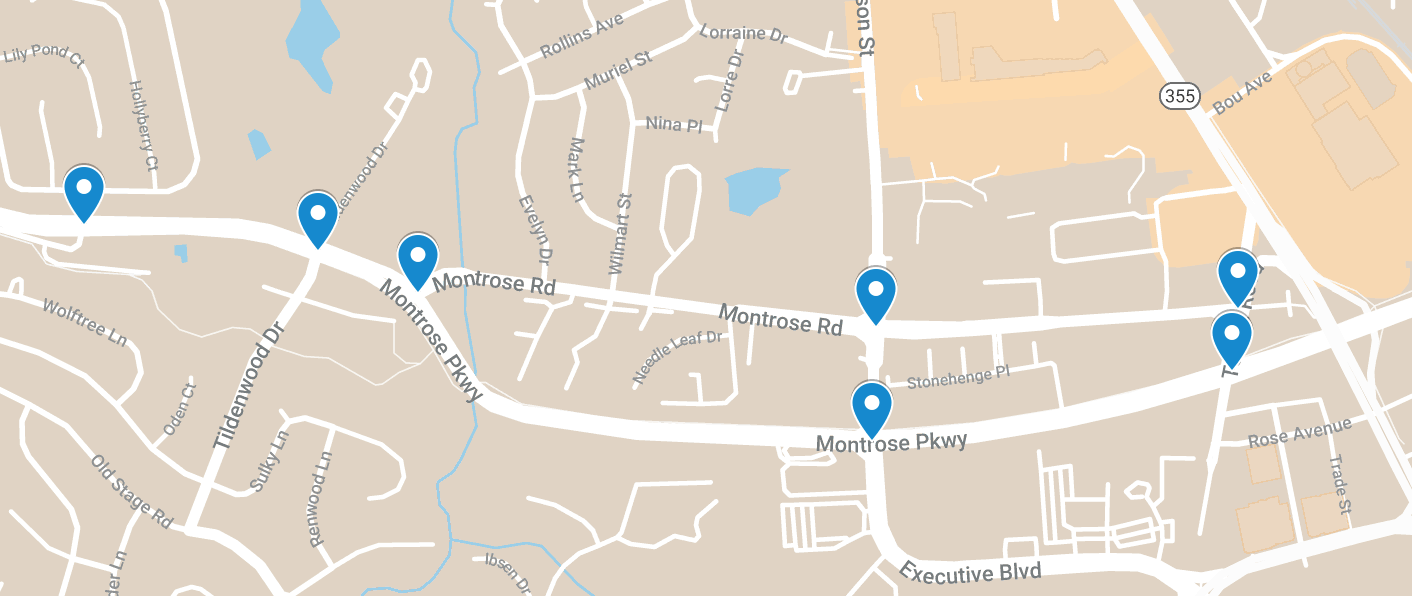}\\[5pt]
\includegraphics[width=.9\columnwidth]{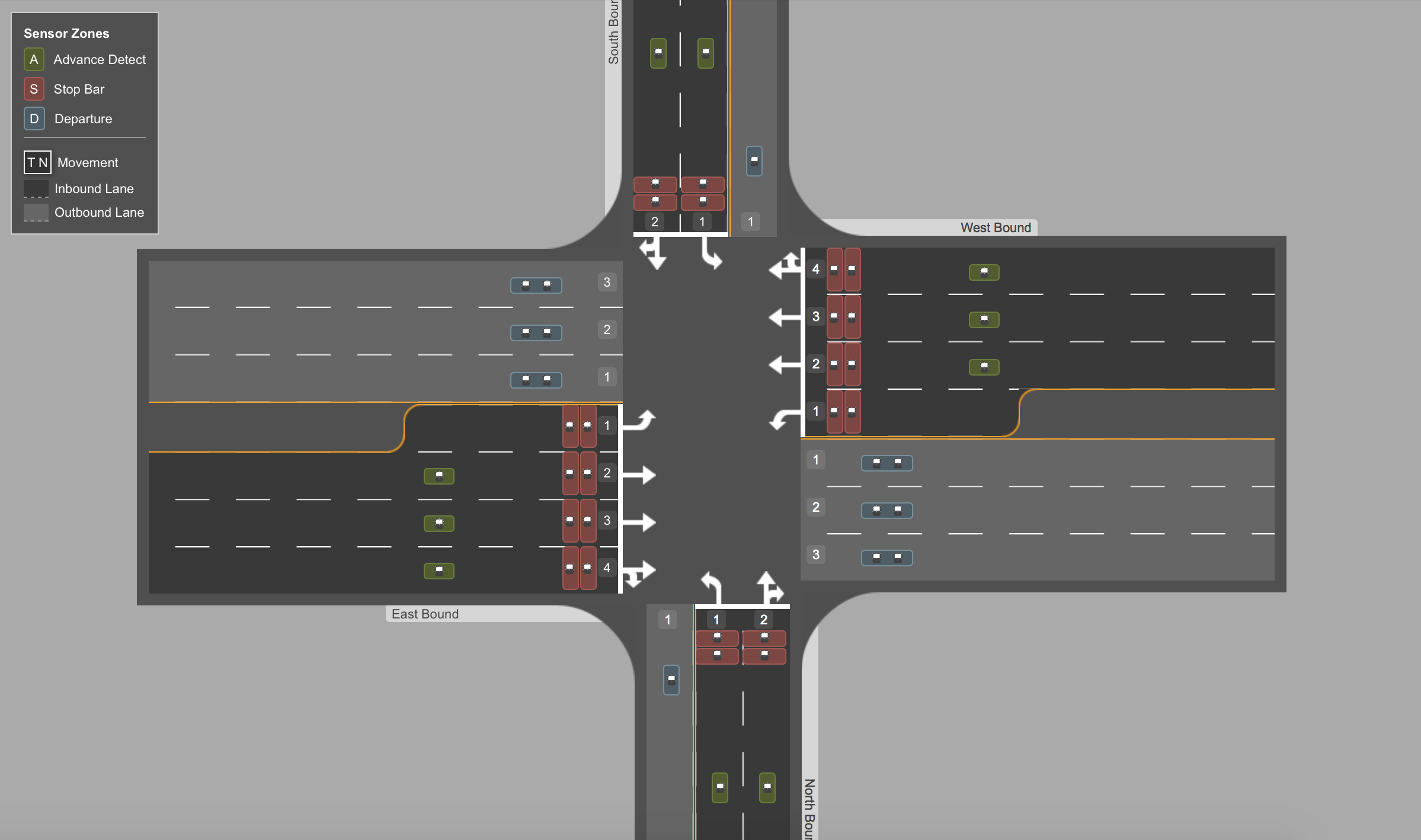}  
  \caption{(Top) A case study traffic network in Maryland consisting of 7 intersections is the basis for the case studies presented below. (Bottom) An example schematic of the sensors at one intersection. Sensors positioned at incoming and outgoing lanes provide high resolution information about vehicle flow rates and queue lengths.}
  \label{fig:dataset}
\end{figure}

As a running case study, this paper considers data collected from a seven intersection network in Maryland. Each intersection is equipped with sensors in incoming and outgoing lanes that detect the presence of vehicles. Vehicle detection events are recorded at the millisecond level, allowing high resolution measurements of traffic flow. In addition, sensors positioned upstream in incoming lanes allow for accurate real-time computation of queue lengths. In the applications below, we will consider both traffic flow and queue lengths as the particular data under consideration. In addition, the traffic light status data is collected, e.g., the currently active phases are known at each time instant. Flow and queue data from a single intersection are used for the Koopman Operator applications in Sections \ref{sec:interpr-modes-sign} and Sections \ref{sec:instability-analysis}. These data are augmented with the traffic signal status information for the application presented in Section \ref{sec:learn-dynam-contr}. Flow data from all intersections in the network are used for the structural inference presented in Section \ref{sec:ident-oper-struct}.


\subsection{Notation}
Using the above case study, in the sequel, we will often consider various data collected over certain lengths of time. In all cases, the data is aggregated into discrete time-steps, usually on the order of several seconds. Let $x_k \in \mathbb{R}^{d_1}$ be the state data variable representing vehicle flows or queues at time $k$ where $d_1$ is the number of sensors or measurement points in the (sub)-network under study. Given $N$ measurement times, $X_1\in \mathbb{R}^{d_1\times (N-1)}$ is a data matrix constructed such that $x_k$ is the $k$-th column of $X_1$ for $k=1,\ldots,N-1$, \emph{i.e.}, The matrix $X_1=\begin{bmatrix}x_1&\cdots x_{N-1}\end{bmatrix}$. $X_2$ is constructed similarly with measurements $k=2,\ldots, N$, \emph{i.e},  $X_2$ is a time-shifted version of $X_1$ containing observations $X_2=\begin{bmatrix}x_2&\cdots&x_N\end{bmatrix}$. 


Similarly, let $u_k \in \{0,1\}^{d_2}$ be the binary input variable representing the signal phases at time $k$ where $d_2$ is the number of possible phases in the (sub)network under consideration. If the $i$-th entry of $u_k$ is 0 ({resp.}, 1), then the $i$-th phase is red (resp., green) at time $k$. Periods when the phase is in yellow or amber is considered part of the green time. $\mathbb{U}$ is the input data matrix whose columns are the input vectors $u_k$ so that $\mathbb{U}=\begin{bmatrix}u_1&\ldots& u_{N-1}\end{bmatrix}$. 


\section{Background Theory and Review}
\label{sec:dmd}
\subsection{Koopman Operator}
Consider a discrete-time system that evolves according to a nonlinear law given by
\begin{equation}
	\begin{split}
	z_{k+1} = \mathcal{T}(z_k)\\
    x_k = f(z_k)
	\end{split}
\end{equation}
\noindent where $z_k\in \mathbb{R}^{n}$ is the state of the system at time $k$, $\mathcal{T}:\mathbb{R}^n\to \mathbb{R}^n$ is a nonlinear vector map describing the evolution of state trajectories, $x_k\in \mathbb{R}^M$ is the measured output, and $f: \mathbb{R}^{n} \to \mathbb{R}^M$ is an observation function, also known as the output function, that maps states to measurements.


Given $\mathcal{T}$, the Koopman operator $\mathcal{K}$ acts on scalar-valued functions of the observation and is given by 
\begin{equation}
	(\mathcal{K} g)(z_k) = g \circ \mathcal{T} = g(z_{k+1}).
\end{equation}
That is, $\mathcal{K}$ is the composition map that composes the function $g$ with the state update map $\mathcal{T}$. Observe that $\mathcal{K}$ is linear, although infinite-dimensional. A function $\phi_j:\mathbb{R}^n\to \mathbb{R}$ along with value $\lambda_j\in \mathcal{C}$ satisfying $\mathcal{K}\phi_j=\lambda_j\phi_j$ is called an \emph{eigenfunction} of $\mathcal{K}$, and $\lambda_j$ is the associated eigenvalue.

Assuming each component of the observation function $f$ lies in the span of the these eigenfunctions, $f$ can be written in terms of the spectral decomposition of $\mathcal{K}$ as 
\begin{equation}
	f(x_k) = (\mathcal{K}^k f)(x_0) = \sum_{j=1}^{\infty} \phi_j (x_k) \psi_j = \sum_{j=1}^{\infty} \lambda^k_j \phi_j(x_0) \psi_j
\end{equation}
\noindent where $\psi_j \in \mathbb{C}^M $ is the \emph{Koopman mode} corresponding to the Koopman eigenvalue $\lambda_j$ belonging to the particular observable $f$. Also, since $x_k$ is vector-valued, we can isolate a particular spatial component and view it in terms of a basis expansion:
\begin{align}
	f(x_k(i)) &= \sum_{j=1}^{\infty} \phi_j (x_k(i)) \psi_j(i), \\
\label{eqn:koopman-expansion}    &= \phi(x_k(i))^T \psi(i), \forall i = \{1,...,M\}.
\end{align}

We see that $\lambda_j$ describes the temporal behaviour of $\psi_j$, where growth or decay is governed by $|\lambda_j|$ and frequency of oscillation by $\angle \lambda_j$. Moreover, from (\ref{eqn:koopman-expansion}), we view $|\psi_j(i)|$ representing the weight of a variable $i$ in a particular mode $j$. On the other hand, $\angle \psi_j(i)$ describes the relative phase of oscillation of the variable.

\subsection{Dynamic Mode Decomposition}
\label{sec:dynam-mode-decomp}
The expansion (\ref{eqn:koopman-expansion}) considers an infinite number of terms because, in general, $\mathcal{K}$ is infinite-dimensional. However, from a computational perspective, it is impractical to consider an infinite expansion of $\mathcal{K}$. Moreover, it has been observed that many high-dimensional systems may be approximated well using a lower-dimensional subspace \cite{maaten_dim_reduction}. Data-driven algorithms seek a finite-dimensional subspace on which to approximate $\mathcal{K}$ \cite{appliedKoopmanism}.  The two main data-driven algorithms that have been connected to approximating $\mathcal{K}$ are dynamic mode decomposition (DMD) \cite{rowley-spectral} and the Arnoldi algorithm \cite{appliedKoopmanism}. In this paper, we use DMD; for a recent book on this method, see \cite{kutz-dmdbook}.

Schmid and Sesterhenn developed DMD in the context of extracting coherent structures in high-dimensional fluid data \cite{Schmid-sesterhenn-dmd, schmid_dmd}. It is formulated as follows: given a sequence of $N$ measurement vectors $x_1,x_2,\ldots,x_N$ with each $x_k\in\mathbb{R}^{M}$, form a pair of time-shifted data matrices $X_1$ and $X_2$ of dimension $\mathbb{R}^{M \times (N-1)}$:
\begin{equation}
\label{eqn:time-shifted}
\begin{split}
X_1 = \begin{bmatrix} 
	  | &  & | \\
      x_1 & \dots & x_{N-1} \\
      | &  & |
	  \end{bmatrix} \space    
X_2 = \begin{bmatrix} 
      | &  & | \\
      x_2 & \dots & x_{N} \\
      | &  & |
	  \end{bmatrix}.
\end{split}
\end{equation}
\noindent where $N$ is the total number of snapshots.

The next step is to find a best fit linear operator $A$ that minimizes the squared Frobenius norm between the two time-shifted matrices in (\ref{eqn:time-shifted}), that is, find $A$ solving
\begin{equation}
\label{eqn:DMD_opt}
\underset{A}{\min}||X_2 - AX_1||_F^2.
\end{equation}
It is known that the minimizing $A \in \mathbb{R}^{M \times M}$ is a finite-dimensional approximation of $\mathcal{K}$ \cite{Tu_dmd}. The minimization (\ref{eqn:DMD_opt}) has a closed form solution that employs the singular value decomposition (SVD) of $X_1$ given by
\begin{equation}
\label{eqn:DMD_sol}
\begin{split}
A &= X_2 X_1^\dagger = X_2 V \Sigma^{-1} U^T \\
\end{split}
\end{equation}
\noindent where $\dagger$ indicates the pseudo-inverse. Since rank$(X_1)= r \leq \min(M, N-1)$, we have $U \in \mathbb{R}^{M\times r}$, $\Sigma \in \mathbb{R}^{r\times r}$, $V \in \mathbb{R}^{(N-1)\times r}$. If $M$ is large, computing the eigen-decomposition of $A$ directly is computationally intensive. To ameliorate this, DMD approximates the eigen-decomposition of a rank-reduced $\tilde{A}\in \mathbb{R}^{r \times r}$ given by
\begin{equation}
\label{eqn:Atilde}
\begin{split}
\tilde{A} &= U^T A U \\
		  &= U^T X_2 V \Sigma^{-1}.
\end{split}
\end{equation}

Further computational gains can be made if there is low-rank structure in the dynamics. In this case, the number of singular values in the SVD of $X_1$ can be truncated to some $\tilde{r} \leq r$ so that
\begin{equation}
\tilde{A} = \tilde{U}^T X_2 \tilde{V} \tilde{\Sigma}^{-1}
\end{equation}
\noindent where $\tilde{U} \in \mathbb{R}^{M\times \tilde{r}}$, $\tilde{\Sigma} \in \mathbb{R}^{\tilde{r}\times \tilde{r}}$ and $\tilde{V} \in \mathbb{R}^{(N-1)\times \tilde{r}}$. The parameter $\tilde{r}$ is a design parameter chosen by, \emph{e.g.}, studying the decay of the singular values.


Given its historical origins, DMD was originally developed for the case when the spatial dimension $M$ of the observation function is much greater than the number of observations $N$. If the converse is true, i.e. $N \gg M$, then the original DMD formulation is inadequate for capturing the nonlinear dynamics fully \cite{Tu_dmd}, \cite{kmd_lowSpatial}.


 One common method to overcome this deficiency is  to enrichen the observables by appending $h$ time-shifted observations to the data matrices \cite{kutz-dmdbook}. In this case, we construct matrices
\begin{equation}
\label{eqn:time-shifted-hankel}
\begin{split}
\tilde{X}_1 = \begin{bmatrix} 
	  | &  & | \\
      x_1 & \dots & x_{N-h} \\
      | &  & | \\
      x_2 & \dots & x_{N-h+1} \\
      | &  & | \\
      \vdots & & \vdots \\
      | &  & | \\
      x_{h} & \dots & x_{N-1} \\
      | &  & | \\
	  \end{bmatrix} \space    
\tilde{X}_2 = \begin{bmatrix} 
      | &  & | \\
      x_2 & \dots & x_{N-h+1} \\
      | &  & | \\
      x_3 & \dots & x_{N-h+2} \\
      | &  & | \\
      \vdots & & \vdots \\
      | &  & | \\
      x_{h+1} & \dots & x_{N} \\
      | &  & | \\
	  \end{bmatrix}
\end{split}
\end{equation}
with $\tilde{X}_1,\tilde{X}_2\in \mathbb{R}^{hM\times (N-h)}$. The parameter $h$ is chosen so that $hM\gg (N-h)$.

\subsection{Adding Control}
\label{sec:adding-control}
Including the effects of  inputs can be critical for systems subject to control, as is the case with traffic networks. In \cite{proctor-dmdc}, DMD is extended to allow for an exogenous control input, and this technique is called dynamic mode decomposition with control (DMDc). DMDc considers the case of a dynamical system with input $u_k\in\mathbb{R}^Q$ acting on the system in time-step $k$. 

DMDc first collects inputs into a matrix $\mathbb{U}$ according to
\begin{equation}
\label{eqn:input}
\begin{split}
\mathbb{U} &= 
\begin{bmatrix} 
| & & | \\
u_1 & \dots & u_{N-1} \\
| & & | \\
\end{bmatrix},
\end{split}
\end{equation}
and then computes $A$ and $B$ so that 
\begin{equation}
\begin{split}
X_2 &= AX_1 + B\mathbb{U} \\
	&= \begin{bmatrix} A & B \end{bmatrix} 
    \begin{bmatrix} X_1 \\ \mathbb{U} \end{bmatrix} \\
    &= G \Omega
\end{split}
\end{equation}
\noindent where $G\in \mathbb{R}^{M\times (M+Q)}$ and $X_1$ and $X_2$ are as previously. As before, if $M<N$, then $\tilde{X}_1$ and $\tilde{X}_2$ can be substituted for $X_1$ and $X_2$. In this case, a matrix $\tilde{\mathbb{U}}\in \mathbb{R}^{hQ\times (N-h)}$ is similarly constructed from $\mathbb{U}$ and used in place of $\mathbb{U}$. Then, $G$ can be computed by
\begin{equation}
\begin{split}
G &= X_2 \Omega^{\dagger} \\
  &= X_2 \hat{V} \hat{\Sigma}^{-1} \hat{U}^T \\ 
\end{split}
\end{equation}
\noindent where $\hat{V}$, $\hat{\Sigma}$ and $\hat{U}$ are from the SVD of $\Omega$ with $\Omega = \hat{U}\hat{\Sigma}\hat{V}^T$.

Since the first $M$ components in the left singular vectors of $\hat{U}$ come from $A$, while the remaining $Q$ components are from $B$, $A$ and $B$ can be determined from $G$ by extracting the appropriate rows in $\hat{U}$ according to
\begin{equation}
\label{eqn:dmdc_learn}
\begin{split}
A &\approx X_2 \hat{V} \hat{\Sigma}^{-1} \hat{U}_1^T \\
B &\approx X_2 \hat{V} \hat{\Sigma}^{-1} \hat{U}_2^T\\
\end{split}
\end{equation}
\noindent where $\hat{U_1} \in \mathbb{R}^{M \times r}$, $\hat{U_2} \in \mathbb{R}^{Q \times r}$ and $\hat{U}=\begin{bmatrix}
\hat{U}_1 \\
\hat{U}_2 \\
\end{bmatrix}$.


\subsection{Overview of Applications}
Here, we provide a brief overview of Koopman Operator theory in various domains, which draw mainly from the ideas of model reduction \cite{appliedKoopmanism} and modal decomposition. The basis for these ideas can be briefly explained as follows: suppose that the first $R \leq M$ Koopman modes capture a sufficient amount of the system's dynamics. Then, a reduced-order model can be approximated as
\begin{equation}
	\label{eqn:koopman-modes-rom}
	f(z_k) = \sum_{j=1}^{R} \lambda_j^k \phi_j (z_0) \psi_j .
\end{equation}

This idea of a constructing a reduced-order model using a partial set of the first $R$ terms leads to the idea of ``separability'' in modes. Then, one could also group terms, for example, based on frequency of oscillation:
\begin{equation}
	\label{eqn:koopman-modes-group}
	f(x_k) = \sum_{j=1}^{R_1} \lambda_j^k \phi_j (x_0) \psi_j + \sum_{j=R_1+1}^{R_2} \lambda_j^k \phi_j (x_0) \psi_j + ...
\end{equation}

\noindent This is called modal decomposition, i.e., partitioning of the system based on the frequency of the modes. Raak et. al \cite{raak-power-partition} test this methodology to partition a power system network. Brunton et. al \cite{brunton_dmdRPCA} separate foreground and background segments in a video, relating slow frequency modes with the background, and fast modes with the foreground.

The notion of growth and decaying modes has also been used in quantitative finance and power systems. In \cite{mann_trading}, the presence of growth or decay modes is used to assess profitable periods for making a trade. In \cite{kmd_power_1}, the number density of unstable Koopman eigenvalues is used to demonstrate short-term stability detection, using data from the 2006 European grid disturbance.

Clustering of modes oscillating at similar frequencies also leads to the idea of coherency analysis. In a power system, maintaining generators at the same synchronous frequency is crucial to the stability of the system. In \cite{susuki-power-coherency}, Koopman modes are used to identify groups of generators oscillating at coherent frequencies. In a disease modeling application, \cite{proctor-disease} use DMD to study how flu and measles spread.


Koopman Operator applications to traffic are limited, with one study involving freeway data \cite{kmd_traffic_1} that compares dominant modes in morning and evening traffic as a means of differentiating the dynamics in those two periods. Drawing from the applications in other domains reviewed here, we demonstrate three novel applications for signalized traffic in this paper. Some preliminary results have been reported in \cite{Ling_KOOPMAN}.


\section{Interpreting Modes for Signal and Phase Timing Estimation}
\label{sec:interpr-modes-sign}
In many cases, it is desirable to estimate traffic signal control parameters such as phase timing directly from traffic flow data. For example, if a legacy traffic control network is retrofitted with wireless sensors, the sensing infrastructure may not be connected to the traffic control hardware and therefore may not have direct access to the status of the traffic controller. In this section, we employ Koopman operator theory to estimate signal and phase timing parameters. These estimates can be used to, \emph{e.g.}, improve energy efficiency of vehicles \cite{lit_optimalControl} and \cite{probabilistic-timing}. 

Koopman Mode analysis is related to use of Koopman Modes to understand the spatio-temporal relationships in a system. In the case of a signalized intersection, the most interpretable spatio-temporal relationship would be the one induced by the phase-split. We apply Koopman Modes for the problem of signal and phase timing estimation, of which motivating examples can be found in \cite{lit_optimalControl} and \cite{probabilistic-timing}. In particular, we show how cycle time, phase sequence and green-splits can be recovered using Koopman Modes. 

\begin{figure}[!ht]
  \centering
  \includegraphics[width=0.5\textwidth]{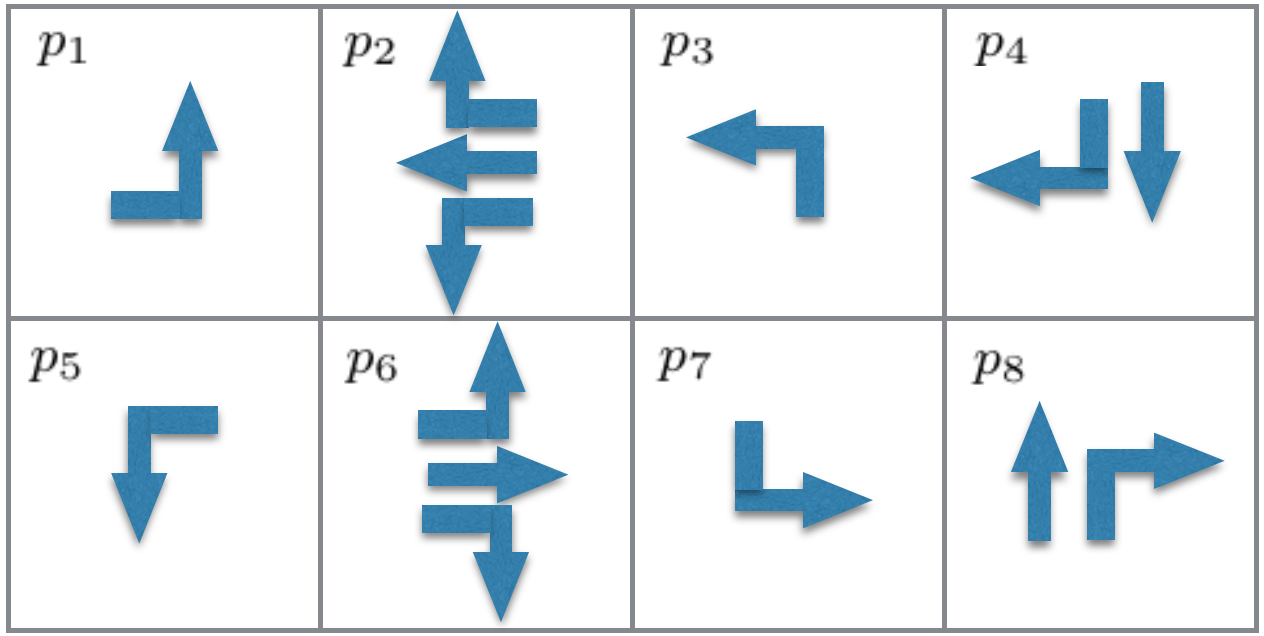}
  \caption{\label{fig:i3-ring} Ring and Barrier Diagram for Montrose Pkwy \& E Jefferson St denoting the phase sequence. Each column of phases is active simultaneously.}
\end{figure}

Our technique is tailored for the standard intersection with 8 phases as shown in Fig. \ref{fig:i3-ring}, however, our basic approach can be generalized to other phase sequencing. We assume the phases in Fig. \ref{fig:i3-ring} are governed by the timing plan in Table \ref{tab:i3-timing}, where each vertical pair of phases have the same green-split. This timing plan is assumed unknown, and we seek to estimate the parameters from traffic flow data. Let $a, b, c$ and $d$ represent the green-splits to be estimated so that $a = p_1 = p_5$s, $b = p_2 = p_6$s, $c = p_3 = p_7$s and $d = p_4 = p_8$s. We first propose a method for estimating the cycle time $C=a+b=c+d$ from a certain dominant Koopman eigenvalue. Then we propose a method for estimating the phase sequence. 


\begin{table}[b]
\caption{Weekday green-splits for Montrose Pkwy \& E Jefferson St (courtesy MCDOT). (Left) 150s cycle from 0600-1000. (Right) 120s cycle from 1000-1500.}
\label{tab:i3-timing}
\centering
\begin{tabular}{|ccccc||ccccc|}
\hline
Phase & $p_1$ & $p_2$ & $p_3$ & $p_4$ & Phase & $p_1$ & $p_2$ & $p_3$ & $p_4$ \\ 
Length (s) & 26 & 53 & 21 & 50 & Length (s) & 16 & 30 & 35 & 39 \\ \hline 
Phase & $p_5$ & $p_6$ & $p_7$ & $p_8$ & Phase & $p_5$ & $p_6$ & $p_7$ & $p_8$  \\
Length (s) & 26 & 53 & 21 & 50 & Length (s) & 16 & 30 & 35 & 39 \\ \hline
\end{tabular}
\end{table}

Let $\lambda_1,\lambda_2,\ldots,\lambda_r$ and $\psi_1,\psi_2,\ldots,\psi_r$ denote a collection of DMD eigenvalues and DMD modes obtained from the traffic flow data as described above.

\paragraph {Estimating Cycle Time}
Let $\lambda_{\theta}$ be the DMD eigenvalue with the largest real part and non-zero imaginary part. Then the cycle time $C$, in seconds, can be estimated according to 
\hfill\break
\begin{equation}
\label{eqn:cycleTimes}
C = \frac{2 \pi \Delta t}{ \mathrm{Im}(\ln(\lambda_{\theta})) }.
\end{equation}

\paragraph {Relative Phase Sequencing and Timing}
Given $\lambda_{\theta}$ used for the cycle time estimation, consider its accompanying mode $\psi_{\theta}$. Recall that $\psi_{\theta}$ is complex-valued and $\{\angle{\psi_{\theta}(1)},...,\angle{ \psi_{\theta}(M)}\}$ provide relative temporal information for each spatial component. We propose using these spatial components to estimate phase sequencing. First, the angles can be converted from radians to seconds by
\hfill\break
\begin{equation}
\label{eqn:phaseLengths1}
\alpha_{\theta}(i) = \text{mod}\left (\frac{\angle{\psi_{\theta}(i)}}{2 \pi} C, C\right).
\end{equation}
for each $i=1,\ldots, M$ where $\text{mod}$ denotes the modulo operator. The modulo operation removes the possibility of $\alpha_{\theta}(i)$ obtaining a negative value. 
Thus, $\alpha_{\theta}(i)$ is a measure of the timing of phase $i$ within the overall cycle length $C$.


\paragraph {Green-Splits}
Given the relative phase timing $\alpha_{\theta}(i)$, we are now in a position to compute the green-splits. From the relative temporal information in $\psi_{\theta}$, the following set of computations are possible: (1) difference between through (T) movements from opposing barriers to compute $b+c$ and $d+a$, (2) difference between left-turn (LT) movements from opposing barriers to compute $a+b$ and $c+d$, and (3) difference between LT and T movements within the same barrier to compute $a$ and $c$. Multiple estimations using different subset combinations can be averaged to determine $a, b, c,$ and $d$.

Since the angles are relative to one another, there is some ambiguity to determine whether smaller angles precede larger angles, or vice versa. To address this, we assume that protected left-turns precede through movements within a leg. Given that the LT angles are larger than the T angles (Fig. \ref{fig:barPhase}), we infer that larger angles precede smaller angles sequentially.

Making use of this information, to compute the difference between movements from opposing barriers, we apply the $\max(\cdot)$ operator on movements from the same barrier. For example, given that both the SBT and NBT movements are green concurrently, as are the EBT and WBT phases, the time difference between T movements can be computed using $\max(\alpha_{\theta}(SBT), \alpha_{\theta}(NBT)) - \max(\alpha_{\theta}(WBT), \alpha_{\theta}(EBT))$. (Note that if the LT angles are smaller than T angles, then we would instead apply the $\min(\cdot)$ operator).

To ensure that the differences are appropriately computed while accounting for $\alpha_{\theta}$ wrapping around with a periodicity of $C$, we check for conditions when one movement is positive while its corresponding movement is negative. For example, if the cycle time is $150$s and $\alpha_{\theta}(SBT)=60$ while $\alpha_{\theta}(NBT)=-55$, the difference $\alpha_{\theta}(SBT) - \alpha_{\theta}(NBT)$ would be incorrectly computed as $115$. Noting that $60$ is equivalent to $-90$, we can adjust $\alpha_{\theta}(SBT)= -\text{mod}(-\alpha_{\theta}(SBT),C)$, and then $\alpha_{\theta}(SBT) - \alpha_{\theta}(NBT)$ would be returned as $35$. The converse could also be true, where $\alpha_{\theta}(NBT)=-55$ is converted to $95$. To disambiguate between the two possibilities, we check if the phases in the opposing barrier are mostly positive or negative. To be precise, denote the two sets $B1$ and $B2$ as:
\begin{equation}
  \begin{split}
    B1 &= \{\alpha_{\theta}(EBT), \alpha_{\theta}(EBLT), \alpha_{\theta}(WBT), \alpha_{\theta}(WBLT) \} \\
    B2 &= \{\alpha_{\theta}(NBT), \alpha_{\theta}(NBLT), \alpha_{\theta}(SBT), \alpha_{\theta}(SBLT)\}. 
  \end{split}
\end{equation}
\noindent Then, if $|\{ B1 | B1>0 \} | > 2$, we would assume that most of the movements in the East-West barrier are positive, and so we would set $\alpha_{\theta}(SBT)=-90$. If the converse is true, then we would set $\alpha_{\theta}(NBT)=95$. We summarize this in Alg. \ref{alg:phases-adjust}.

\begin{algorithm}[!ht]
\caption{Adjust Angles} \label{alg:phases-adjust}
\begin{algorithmic}[1]
\scriptsize
    \IF{$\alpha_{\theta}(SBT) - \alpha_{\theta}(NBT) \geq C/2$}
        \IF{$| \{ B1 | B1 > 0 \} | > 2$}
            \STATE $\alpha_{\theta}(SBT) = -mod(-\alpha_{\theta}(SBT), C)$
        \ELSE
            \STATE $\alpha_{\theta}(NBT) = mod(\alpha_{\theta}(NBT), C)$
        \ENDIF    
    \ELSIF{$\alpha_{\theta}(SBT) - \alpha_{\theta}(NBT) \leq -C/2$}
        \IF{$| \{ B1 | B1 > 0 \} | > 2$}
            \STATE $\alpha_{\theta}(NBT) = -mod(\alpha_{\theta}(-NBT), C)$
        \ELSE
            \STATE $\alpha_{\theta}(SBT) = mod(\alpha_{\theta}(SBT), C)$
        \ENDIF 
    \ENDIF
    \STATE Repeat steps 1-13 for SBLT and NBLT
    \STATE Repeat steps 1-13 for WBT and EBT
    \STATE Repeat steps 1-13 for WBLT and EBLT
\end{algorithmic}
\end{algorithm}

\subsection{Case Study}
\begin{figure}[!t]
  \centering
  \includegraphics[width=0.5\textwidth]{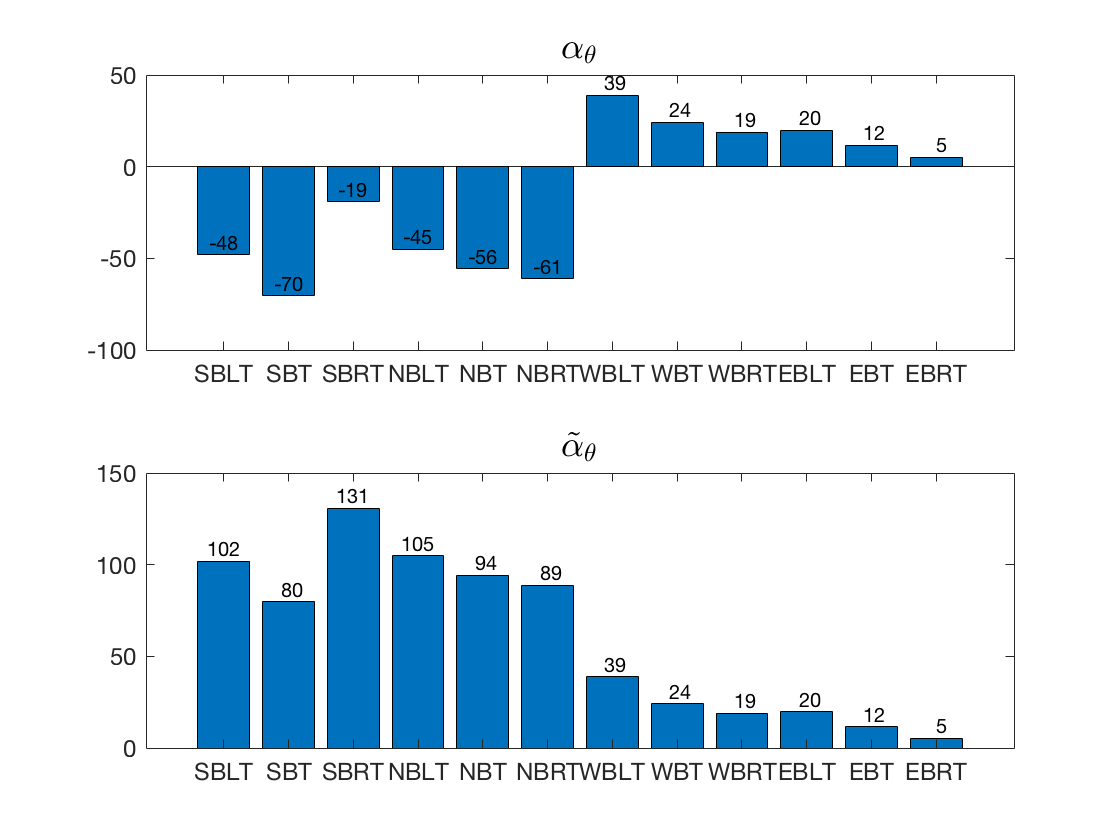}
  \caption{\label{fig:barPhase} Phase-splits, in seconds, of the 150s mode for vehicle flow from 0900-0930 on 14 February 2017 computed using angles obtained from the corresponding Koopman Mode.}
\end{figure}

\begin{figure}[!t]
    \centering
    \includegraphics[width=1.\linewidth]{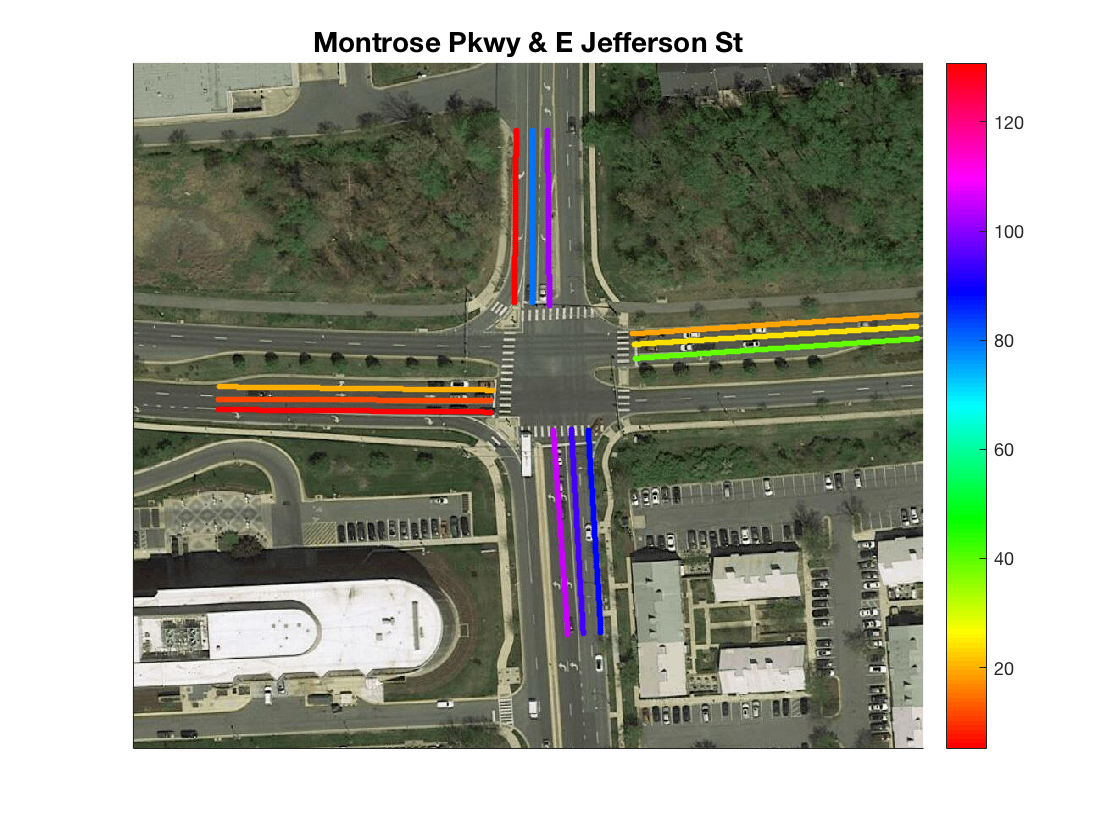}
  \caption{\label{fig:Sequence}Phase-splits, in seconds, corresponding to the 150s mode, when running DMD using samples from 0900-0930. Lines indicate direction, e.g. the lines farthest left are the EB movements, with the top line as the EB-LT movement, while the bottom line on the farthest right is the WB-LT movement.}
\end{figure}

We demonstrate the above estimation technique for cycle time, phase sequence, and green split timing on a case study using vehicle flows at the Montrose Pkwy \& E Jefferson St Intersection. Based on the timing plan for this intersection, on weekdays, it is known that there are two cycle times: 120s (for the hours 1000-1500 and 1900-0000) and 150s (for the hours 0600-1000 and 1500-1900), while during weekends, there is only a 120s cycle. Our goal is to estimate these cycle times, and the accomponing phase sequences and green splits, using only traffic flow data. To achieve this, we perform DMD on vehicle flows binned at 10s from the Montrose Pkwy \& E Jefferson St. intersection. 

First, we run an experiment on a test week of 12-18 February 2018, making hourly predictions of the cycle time rounded to the nearest integer, from 0600-2200 using $N=360$ (Table \ref{tab:fft-exp2}). The cycle time is correctly estimated within a tolerance of $\pm3$ seconds with high accuracy on all the days. Here, we find that low and inconsistent vehicle flow affects the accuracy of the prediction, where the incorrect cycles beyond a tolerance of $\pm3$ seconds frequently occur during the 0600-0700 time. Also, traffic flow on weekends tend to be inconsistent and of lower volume, and so they have a lower number of predicted cycles returned to the exact value compared to weekdays.

Next, we interpret phase sequence from $\tilde{\alpha}_{\theta}$. We plot $\alpha_{\theta}$ and $\tilde{\alpha}_{\theta}$ in Fig. \ref{fig:barPhase}. We additionally overlay $\tilde{\alpha}_{\theta}$ in a Google Map of the intersection (Fig. \ref{fig:Sequence}). Referring to this figure, the phase sequence can be inferred as: (SBLT and NBLT) $\rightarrow$ (SBT and NBT) $\rightarrow$ WBLT $\rightarrow$ WBT $\rightarrow$ EBLT $\rightarrow$ EBT, which matches the Ring and Barrier Diagram in Fig. \ref{fig:i3-ring}. We note that using this method, turn-movements associated with multiple phases (e.g. LT movements in permissive-protected mode) may appear immediately before, during, or immediately after its associated Through movement in the sequence, depending on the arrival of vehicles. For example, if no vehicles are present during the protected mode of the LT and arrive only during the permissive phase, then the LT and T phases would have approximately the same color. Conversely, if vehicles arrive during the protected mode of the LT, then the LT phase would have a color that is sequentially before the T phase. For this time period we observe the latter case, as all the LT's precede the T's. RT's are ignored in the analysis, since they are typically allowed on red lights.

We next estimate the green-split for the 1000-1500 time (120s cycle), applying Alg. \ref{alg:phases} using $N=360$, and average the results made from applying a sliding window of 10s in each test interval. Also, we choose $h$ to be an integer multiple of C with a minimum of 2 cycles, such that $h=24$ for the 120s cycle. We omit each iteration where the estimated cycle time is not exactly 120s. A table of the results is provided in Tab. \ref{tab:green-results-120s}, which show reasonable accuracy. We highlight that despite using vehicle flows binned at 10s, the algorithm is able to recover the green-splits at second-level precision.

\begin{algorithm}[!ht]
\caption{Estimate Green-Split Lengths} \label{alg:phases}
\begin{algorithmic}[1]
\scriptsize
    \STATE \textbf{LT's from opposing barriers:} \\
    \STATE $\kappa = max(\alpha_{\theta}(SBLT), \alpha_{\theta}(NBLT))- max(\alpha_{\theta}(EBLT), \alpha_{\theta}(WBLT))$ \\
    \STATE $c+d = mod(\kappa, C)$ \\
    \STATE $a+b = C - (c+d)$
    \STATE \textbf{T's from opposing barriers:} \\
    \STATE $\kappa = max(\alpha_{\theta}(SBT), \alpha_{\theta}(NBT)) - max(\alpha_{\theta}(EBT), \alpha_{\theta}(WBT))$ \\
    \STATE $b+c = mod(\kappa, C)$ \\
    \STATE $d+a = C - (b+c)$
    \STATE \textbf{LT and T within E-W barrier:}
    \STATE $\tilde{a} = max(\alpha_{\theta}(EBLT), \alpha_{\theta}(WBLT)) - min(\alpha_{\theta}(EBT), \alpha_{\theta}(WBT))$ 
    \STATE $a = mod(\tilde{a}, C)$
    \STATE \textbf{LT and T within N-S barrier:}
    \STATE $\tilde{c} = max(\alpha_{\theta}(SBLT), \alpha_{\theta}(NBLT)) - min(\alpha_{\theta}(SBT), \alpha_{\theta}(NBT))$ 
    \STATE $c = mod(\tilde{c}, C)$
\end{algorithmic}
\end{algorithm}

To summarize this section, we have demonstrated Koopman Modes for analyzing spatio-temporal relationships at a traffic intersection, in particular, for signal and phase timing recovery. 

\begin{table}[t]
\caption{Accuracy for recovering cycle time in hourly intervals from 0600-2200 on a test week at the Montrose Pkwy \& E Jefferson St Intersection, at different tolerance intervals. All cycle time estimates are rounded to the nearest integer.}
\label{tab:fft-exp2}
\centering
\begin{tabular}{|c|c|c|c|}
\hline
Day & Exact & Tol. +/- 1s & Tol. +/- 3s \\
& (rounded to nearest integer) & & \\ \hline \hline
12-Feb-2017 & 9/17  & 14/17 & 15/17 \\ 
13-Feb-2017 & 13/17 & 15/17 & 16/17 \\
14-Feb-2017 & 14/17 & 15/17 & 16/17 \\
15-Feb-2017 & 13/17 & 16/17 & 17/17 \\
16-Feb-2017 & 10/17 & 15/17 & 15/17 \\
17-Feb-2017 & 12/17 & 16/17 & 16/17 \\
18-Feb-2017 & 6/17 & 15/17  & 17/17 \\
\hline
\end{tabular}
\end{table}

\begin{table}[t]
\caption{Error for estimated green lengths over several time intervals for 120s cycle.}
\label{tab:green-results-120s}
\centering
\begin{tabular}{|c|c|c|c|c|c|}
\hline
Error  & Time & $a$ & $b$ & $c$ & $d$ \\
& Interval & & & &  \\ \hline \hline
sec & 1000-1100   & -7         & -1         & 6          & 2     \\ 
\% &            & 43.8       & 3.3        & 17.1       & 5.1   \\ \hline
sec & 1100-1200   & -2         & -1         & 6          & -3    \\ 
\% &            & 12.5       & 3.3        & 17.1       & 7.7    \\ \hline
sec & 1200-1300   & -3         & 0          & 5          & -2     \\ 
\% &            & 18.8       & 0          & 14.3       & 5.1   \\ \hline
sec & 1300-1400   & 0          & -1         & 6          & -5    \\ 
\% &            & 0          & 3.3        & 17.1       & 12.8   \\ \hline
sec & 1000-1400   & -2         & 0          & 5          & -3     \\ 
\% &            & 12.5       & 0          & 14.3       & 7.7   \\ \hline
\end{tabular}
\end{table}


\section{Instability Analysis}
\label{sec:instability-analysis}
If traffic at an intersection becomes too congested, the queues that accumulate during the red phase may not fully dissipate during the next green phase, leading to unstable queue growth. The spectral decomposition of the Koopman Operator allows for identifying and analyzing such instability. In particular, knowledge of growth modes from traffic queue-dynamics can supplement knowledge of existing queue lengths for anticipating congestion. In this way, while presence of high queue volumes is indicative of a traffic incident, detecting sustained queue growth would serve as an early warning. We propose using DMD to learn the queue dynamics on a rolling window of queue measurements. 


Let the snapshot vector $x_k$ be the queue length at a particular leg, with each increment of time-step $k$ corresponding to $\Delta t = 10$ s. Form the data matrices $\tilde{X}_1$ and $\tilde{X}_2$ according to (\ref{eqn:time-shifted-hankel}), using $h=10$. Then, run DMD in a rolling window fashion with a stride of $k$, looking back at the $N$ most recent samples. Apply a rank truncation of $\tilde{r}=10$. Also, let $\lambda_{1,k}$ be the largest eigenvalue of $A$ at time $k$, $\Upsilon_k$ the count of consecutive $|\lambda_1|>1$ at time $k$, and $\tilde{N}$ the total number of samples considered in the day. We then propose a thresholding scheme on $\Upsilon_k$ to indicate detected queue instability. Alg. \ref{alg:instability} summarizes this process.

\begin{algorithm}[!t]
  \caption{Instability Detection Algorithm} \label{alg:instability}
  \begin{algorithmic}[1]
      \scriptsize
      \STATE Initialize $ctr = 0$ and $N=180$.
      \FOR{$k = N: \tilde{N}$}
      	\STATE Run DMD to learn $A$ using observations from $k-N$ to $k$.
	  \IF{$|\lambda_{1,k}|>1$}
      	\STATE $ctr = ctr + 1$
      \ELSE
      	\STATE $ctr = 0$
      \ENDIF
      \STATE $\Upsilon_k=ctr$.
      \IF{$\Upsilon_k>\epsilon$}
      	\STATE Flag an abnormal event condition.
      \ENDIF
       \ENDFOR
  \end{algorithmic}
\end{algorithm}

\begin{figure*}[!t]
\subfloat[West-Bound Leg, Normal Day]{\includegraphics[width=3.7in]{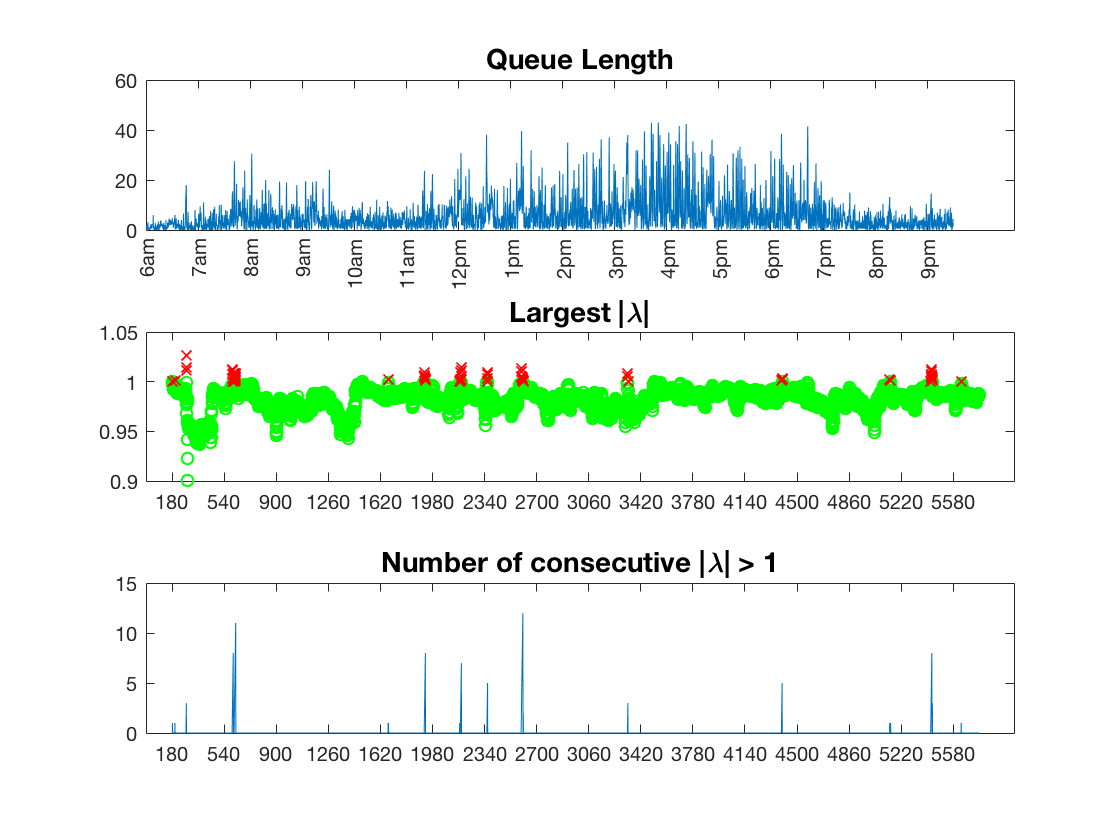}
\label{fig:eig-normal}}
\hfil
\subfloat[West-Bound Leg, Accident Day]{\includegraphics[width=3.7in]{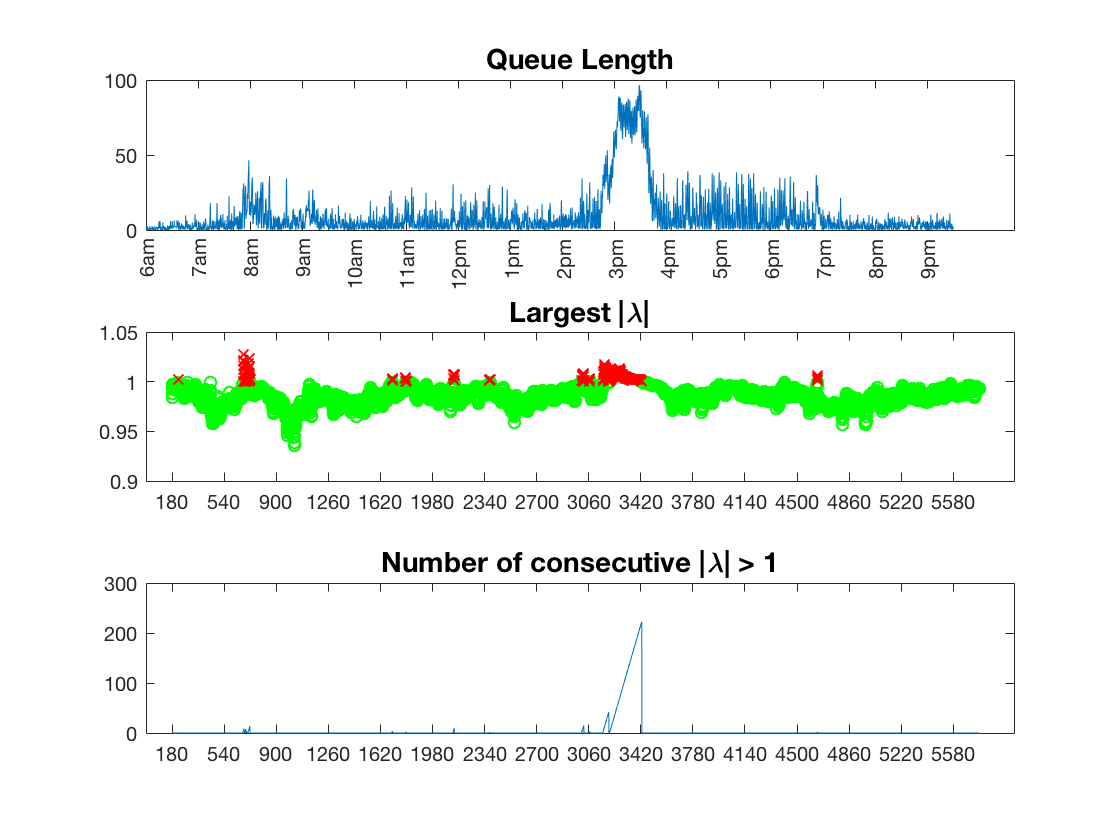}
\label{fig:eig-accident}}
\caption{(Top Subplots) Queue Lengths. (Middle Subplots) $|\lambda_{1,k}|$ (unstable $\lambda_{1,k}$ are shown in red, while stable $\lambda_{1,k}$ are green). (Bottom Subplots) $\Upsilon_k$, the number of consecutive unstable eigenvalues. Each data-point in the Middle and Bottom subplots corresponds to the dynamics of the previous 180 data samples. Notice that the number of consecutive unstable eigenvalues is much larger during the accident.}
\label{fig:unstable_rolling}
\end{figure*}

To demonstrate our approach, we consider data from Friday, 10 February 2017. At approximately 2.47pm on this day, there was a reported accident at the Tildenwood-Montrose Road Intersection which affected the WB and EB movements. We provide a plot of the queues at the WB movement, shown in the top subplot in Fig. \ref{fig:eig-accident}. (The accident data is obtained from the Maryland Open Data Portal \cite{maryland}). Notice the prolonged spike from the time of the accident until about 4pm, which is more severe compared to the peak hour in the morning. For comparison, we also provide queue lengths at the same movement on a normal Friday, 17 February 2017 in the top subplot in Fig. \ref{fig:eig-normal}.

Using Alg. \ref{alg:instability}, we compare the number of unstable eigenvalues for the normal day and the accident day. We choose $N=180$, which corresponds to 30 minutes of data binned at 10s intervals. Choice of $N$ is related to how long of a window to model the dynamics. Since accidents may have sustained growth compared to peak hour growth, a fairly long window such as 30 minutes would be appropriate.

In both Fig. \ref{fig:eig-normal} and Fig. \ref{fig:eig-accident}, the top subplot shows queue lengths at the WB leg, the middle subplot shows $|\lambda_1|$ obtained at each time-step, where unstable $\lambda_1$ are shown in red while stable $\lambda_1$ are green, and the bottom subplot shows $\Upsilon_k$. Notice that the number of consecutive unstable eigenvalues is at least five times larger during the accident compared to normal peak hours. At $k=3420$, $\Upsilon_{3420} \approx 200$ in Fig. \ref{fig:eig-accident}. To provide intuition on how to interpret $\Upsilon_k$, since a stride of $10s$ is used, this corresponds to about 33 minutes of sustained queue growth as of time-step $k=3420$.

\begin{figure}[!t]
\centering
\includegraphics[width=3.75in]{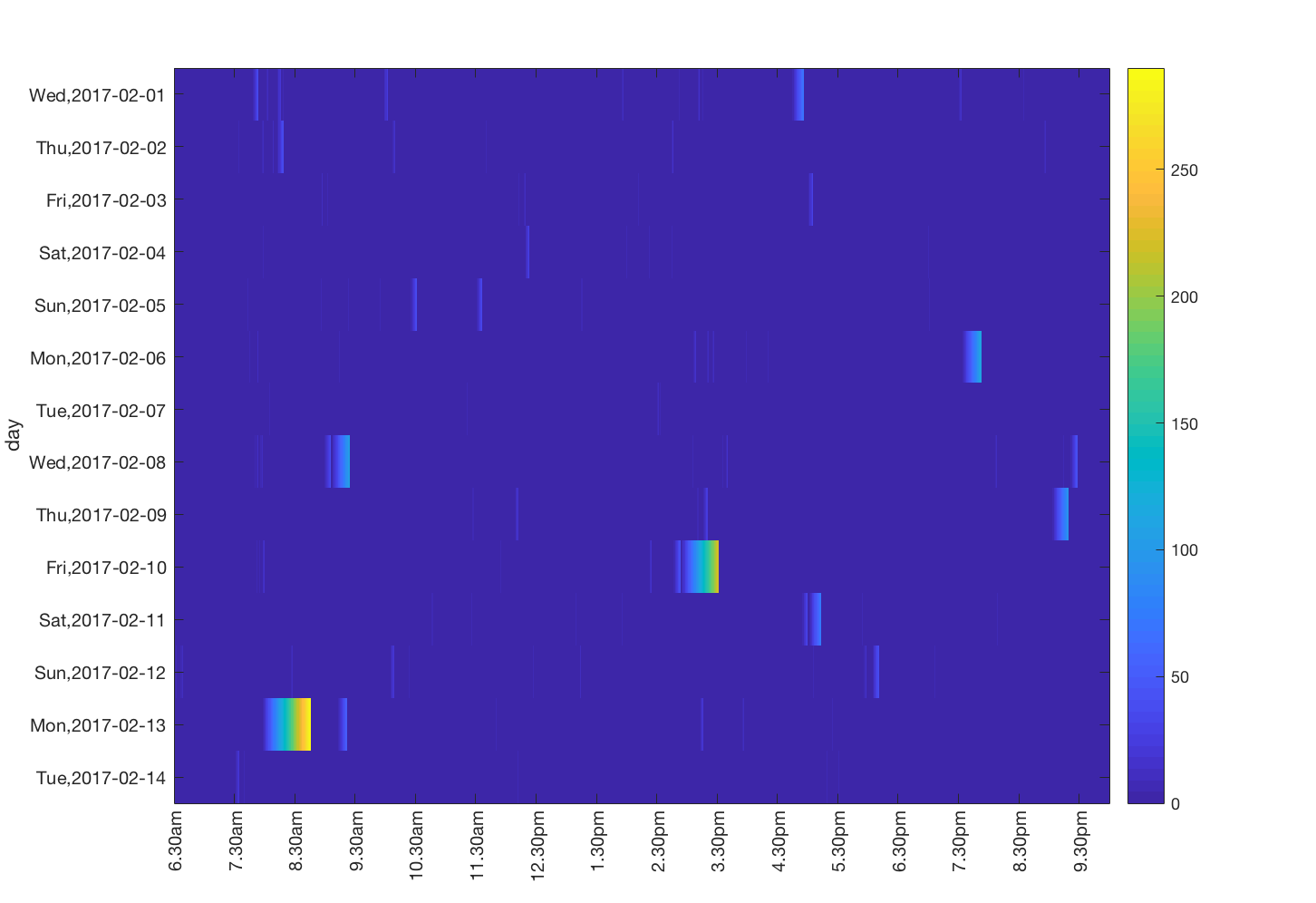}
\caption{Consecutive unstable eigenvalues for the WB Leg at Montrose Rd \& Tildenwood intersection over 2 Weeks.}
\label{fig:2Weeks-Sequence}
\end{figure}

The algorithm can also be used in a post-analysis study to identify recurring and/or abnormal periods of congestion. To demonstrate this, we provide a visualization of the algorithm run on the same intersection over a period of 2 weeks (Fig. \ref{fig:2Weeks-Sequence}). From the figure, two periods of increased growth can be identified as shown by the yellow regions; 10 February 2017 at around 2.40pm, and 13 February 2017 at around 7.40am. The increased growth on 10 February as described earlier is related to an accident. We are unable to find a recorded incident for congestion on the second day, but postulate that it could be due to lane closures or simply an unusual congestion from which no record can be easily found online. Nevertheless, a record of current traffic incidents could be mined from social websites such as Twitter \cite{lit_twitterIEEE}, \cite{lit_twitterTRB} and be compared against identified growth periods in the data.


\section{Learning Dynamics for Control}
\label{sec:learn-dynam-contr}
\begin{figure*}[!t]
\centering 
\subfloat[0\% green time for EB-WB phases]{\includegraphics[width=3.in]{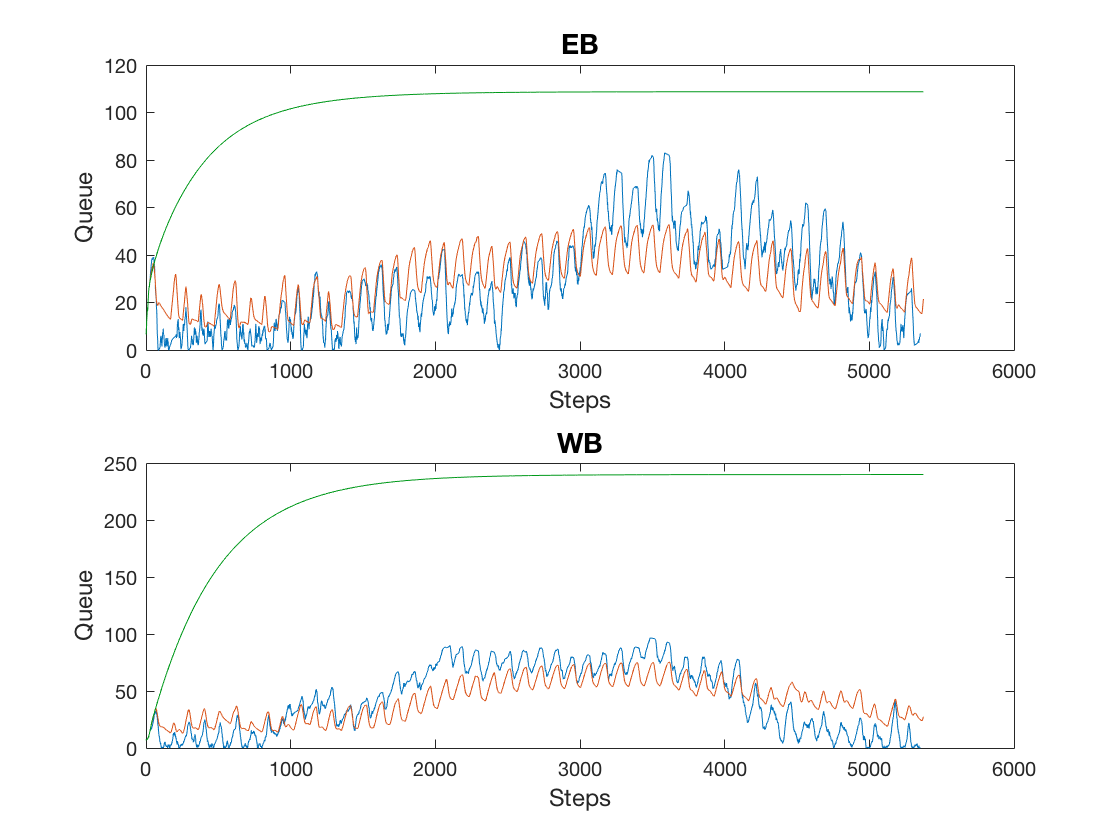}%
\label{fig:Greens0}} 
\hfil
\subfloat[20\% green time for EB-WB phases]{\includegraphics[width=3.in]{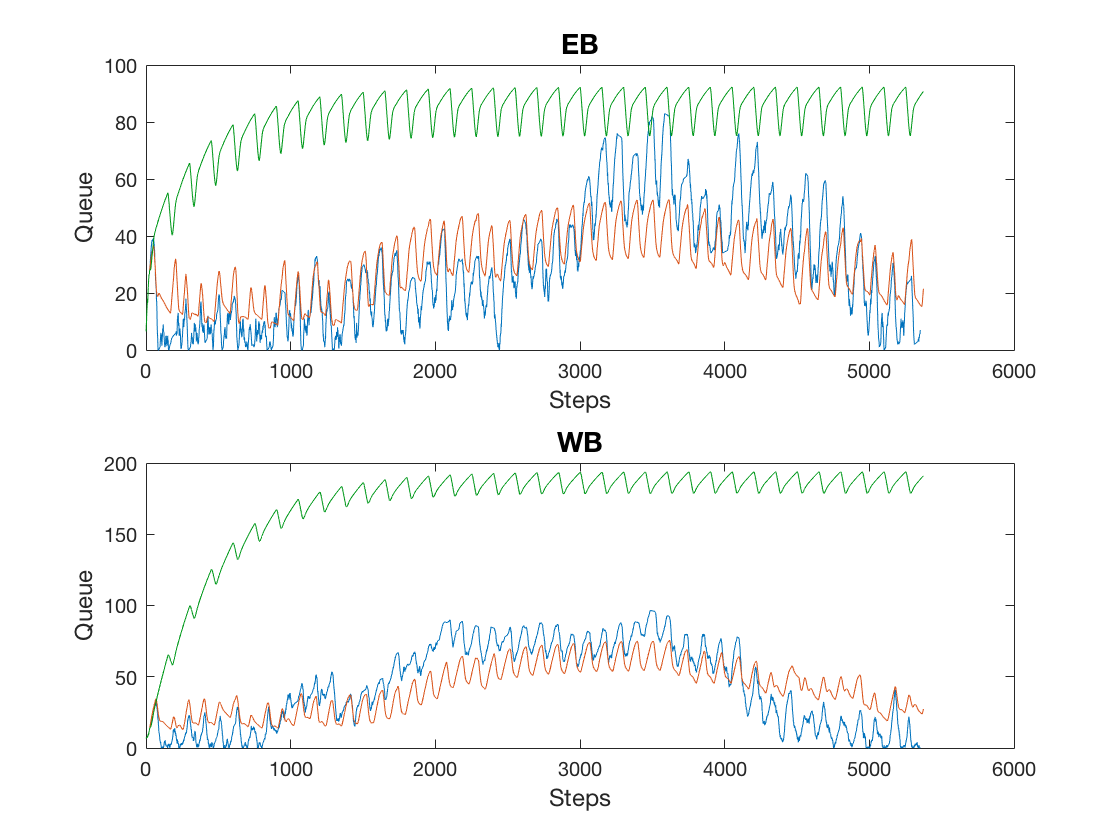}%
\label{fig:Greens20}} 
\hfil
\subfloat[80\% green time for EB-WB phases]{\includegraphics[width=3.in]{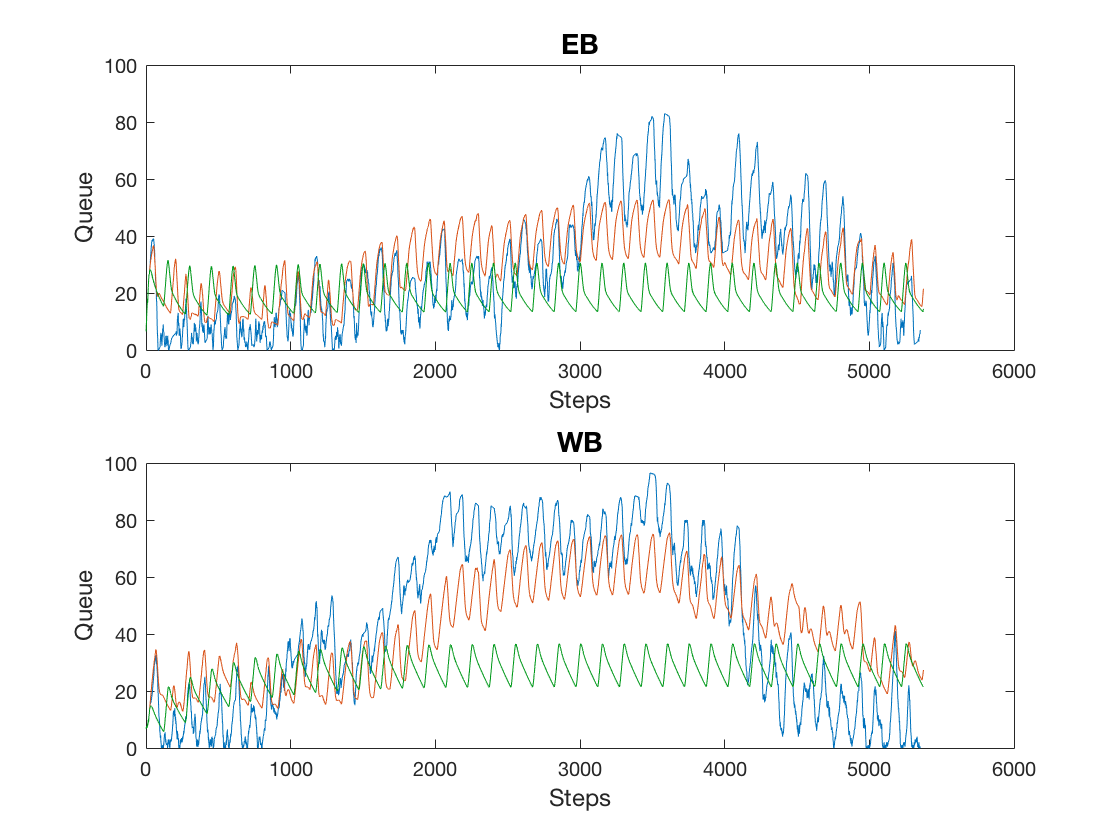}%
\label{fig:Greens80}} 
\hfil
\subfloat[100\% green time for EB-WB phases]{\includegraphics[width=3.in]{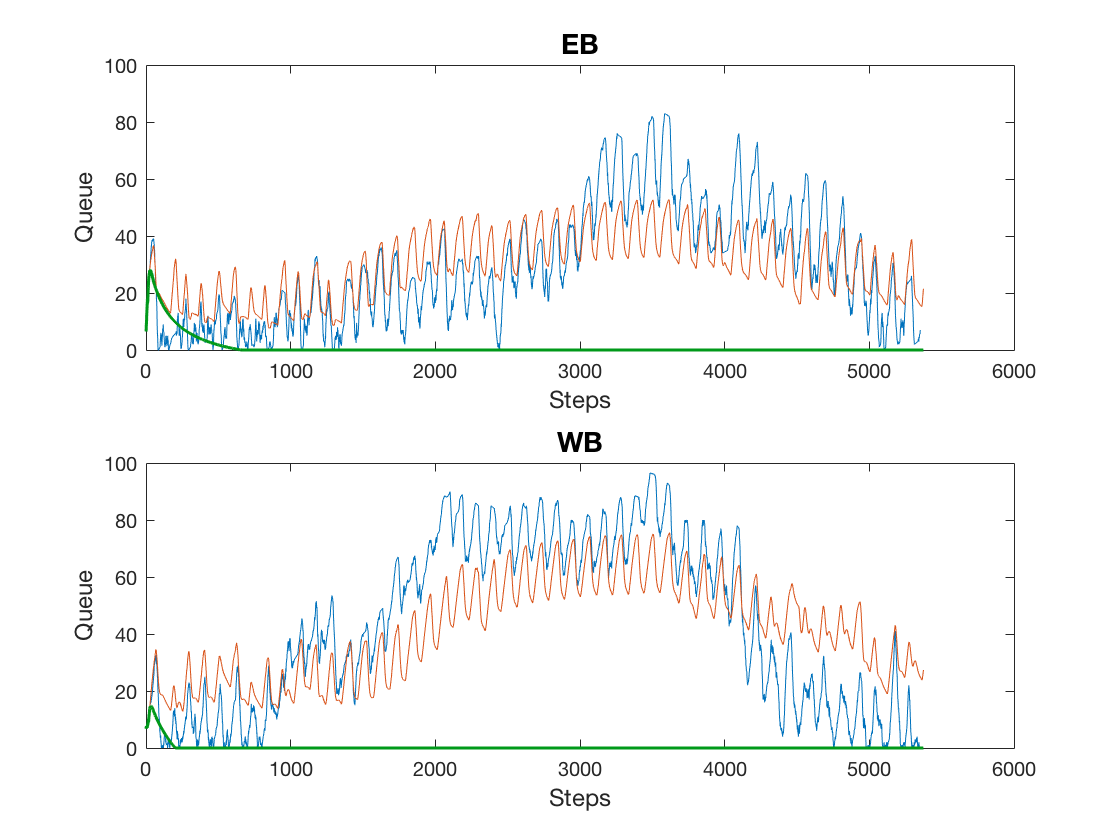}%
\label{fig:Greens100}} 
\caption{The effect of varying green percentages for the EB-WB phases at Montrose Rd \& Tildenwood intersection. When 0\% green time is allocated, the queue saturates at the queue capacity, while when 100\% green time is allocated, the queue decays quickly to zero. (Green) Predicted queues under modified green-split. (Red) Predicted queues under original green-split. (Blue) Original queues.}
\label{fig:dmdc-sp-switch}
\end{figure*}
In the previous section, Koopman operator theory was used to identify unstable queue growth, and it was observed that this approach can identify early stages of instability. This allows traffic operators to be alerted early with the hope that appropriate adjustments can be made to signal timing in order to mitigate the instability. It is also possible to use Koopman operator theory to approximate the queueing process under alternative control schemes. In this section, we propose using DMDc, introduced in \cite{proctor-dmdc} and discussed in Section \ref{sec:adding-control}, to learn an approximate model of the queue dynamics using signal phases as exogenous input. Thus, in addition to assuming queue measurements, we assume access to the traffic signal state at each time instant.
In \cite{Ling_KOOPMAN}, we used DMDc with the aim of evaluating anticipated queue lengths under a different timing scheme than the pre-programed split during congestion. In this paper, we extend the analysis to study the effect of the queues under different green time percentages.

Define the input vector $u_k$ as the signal phase at each turn movement at the affected intersection i.e.,
\begin{equation*}
u_k = \begin{bmatrix}
	u_k(1) \\ u_k(2) \\ \vdots \\ u_k(12)
	\end{bmatrix}
\end{equation*}
where \begin{equation*}
		u_k(i) = 
			\begin{cases}
      			0, & \text{if } u_k(i) = \text{red}   \\
      			1, & \text{if } u_k(i) = \text{green or yellow}.
			\end{cases}
	  \end{equation*}

Let the state vector $x_k$ at time $k$ consist of aggregated queue lengths of each leg at the intersection so that $x_k \in \mathbb{R}^4$. Then, form the data matrices $\tilde{X_1}$, $\tilde{X_2}$ and $\tilde{\mathbb{U}}$ respectively, as in \eqref{eqn:time-shifted-hankel} and \eqref{eqn:input}, using $h = 12$. We use a rank truncation $\tilde{r}$ based on the number of singular values above $1e^{-10}$.

First, we estimate the $A$ and $B$ matrices using (\ref{eqn:dmdc_learn}) where $u_k$ consists of the original signal phase. Here, $x_1$ is the vector of queues at 2.30pm, $x_N$ is the vector of queues at 4pm, and similarly for $u_1$ and $u_N$. To validate that $A$ and $B$ encode an appropriate model of the queue dynamics, we reconstruct the queues in multi-step fashion using $x_1$ as the initial condition and the original $\{u_1,...,u_N\}$ as in Alg. \ref{alg:reconstruction}. The reconstructed queues for the EB and WB legs (shown in red) in Fig. \ref{fig:dmdc-sp-switch} show a good approximation of the dynamics of the original queues (shown in blue).

\begin{algorithm}
  \caption{Multi-Step Reconstruction} \label{alg:reconstruction}
  \begin{algorithmic}[1]
      \scriptsize
      \STATE Inputs: $A$, $B$, $\{u_1,...,u_N\}$ and $x_1$.
      \STATE Returns: $\{x_2,...,x_N\}$
      \FOR{$k = 1: N$}
      \STATE $x_{k+1} = Ax_k + Bu_k$
      \ENDFOR
  \end{algorithmic}
\end{algorithm} 

In preliminary results reported in \cite{Ling_KOOPMAN}, we have investigated the effect of modified inputs for two cases: (i) longer green-splits for the WB-EB phases, and (ii) shorter green-splits for the WB-EB phases, using signal phases from the recorded actuations. Our main findings are that under the longer green-split, the queues for the WB-EB legs are less severe while under the shorter green-split, there is persisted high volume of queues at the WB-EB legs.

In this paper, we investigate the effect of the signal phases on predicting queue dynamics. Examining the structure of the $B$ matrix, we observe that $B$ encodes the effect of the light status on the original queue dynamics. We note that there are significant non-zero entries are in the last $d_1$ rows of $B$. Let $b_n^T$ denote one of the last $d_1$ rows of B, while $u_n$ denote one of the last $d_1$ columns of $\mathbb{U}$. Then, the inner product $b_n^T u_n$ expresses the contribution of the signal phase to the queue $x_n$, or the $n^{th}$ row in $X_2$.

If $u_n(i)=0$ (the phase is red), then it does not matter what the value of $b_n(i)$ is. However, if $u_n(i)=1$ (the phase is green or yellow), there are 3 possible cases: (i) $b_n(i) > 0$, where the green light contributes to an increase in queue length at leg $i$, (ii) $b_n(i)=0$, where there is no net contribution due to the light being green, and (iii) $b_n(i) < 0$, where the green light contributes a reduction in queues at leg $i$. By increasing the effective green time of a leg, there is a longer sequence of $u_n(i)=1$, thus providing overall an increased contribution of $b_n(i) < 0$, thus effectively reducing the queue length.

We demonstrate reconstructed queue lengths for the EB-WB leg under different green-split percentages; $0\%$, $20\%$, $80\%$, and $100\%$ in Fig. \ref{fig:dmdc-sp-switch}. For these experiments, we create the input phase matrix synthetically, instead of using a sample of the signal phases from the data. When the EB-WB legs receive $0\%$ green time, the queues grow until saturating at what could be the learned queue capacity. Conversely, under $100\%$ green time, the queues decay quickly to $0$. The queues receiving $20\%$ and $80\%$ green time respectively show behavior in between these two extremes.

We remark that although the reconstructed queues may not be exact, they provide qualitative insight on how the queue dynamics change under different timing schemes. Also, while the input signal phases are at turn-movement resolution, the queue lengths we have utilized are at leg-resolution. 
\section{Identifying Network Structure}
\label{sec:ident-oper-struct}
In the previous sections we have focused on applications of DMD for single intersections. 
In this section we show that DMD can be used to infer geometric structural relationships in a network of signalized intersections. We achieve this by investigating the $A$ matrix of DMD satisfying \eqref{eqn:DMD_opt}. Here, we also employ the time-shifting augmentation technique suggested in  \eqref{eqn:time-shifted-hankel}. Inspired by the special structure induced by this approach, we further compare this approach to vector autoregression for prediction. Our results indicate that DMD has certain advantages over vector autoregression.

\subsection{DMD Matrix Structure}
As described in Section \ref{sec:dynam-mode-decomp}, if the spatial dimension $M$ is much lower than the number of observations $N$, then it is often beneficial to increase the rank of the data matrix by augmenting the system with delayed as shown in \eqref{eqn:time-shifted-hankel}. In general, we may
increase the number $h$ of delay coordinates to allow a larger operator matrix $A$ and more Koopman modes. As an example, consider queue length data from the Montrose Rd \& E Jefferson St. intersection measured at each second, and let $x_k$ consist of aggregated queue lengths of each leg at the intersection at time $k$ so that $x_k \in \mathbb{R}^4$. We choose $N=1200$ and $h=10$, meaning that we use $1200$ consecutive observations (20 minutes) to form the data matrices $\tilde{X_1}$, $\tilde{X_2}$, as in \eqref{eqn:time-shifted-hankel}, after including nine additional time-shifted block-rows. We do not discard any modes using a rank truncation.

\begin{figure}[t]
	\centering
	\includegraphics[width=0.8\linewidth]{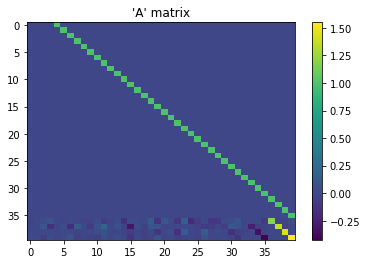}
	\caption{\label{fig:a} A matrix structure from DMD with $10$ time-shifted data. Each data point has dimension 4 so the A matrix has dimension $40 \times 40$.}
\end{figure}

Figure \ref{fig:a} shows the heat map of the $A$ matrix from \eqref{eqn:DMD_sol}. Observe that $A$ tends to have the block structure 
\begin{equation}
\label{eqn:structure}
A = \begin{pmatrix}
0 & I & 0 & \cdots & 0 \\
0 & 0 & I & \cdots & 0 \\
\vdots & \vdots & \vdots & \ddots & \vdots\\
0 & 0 & 0 & \cdots & I \\
A_1 & A_2 & A_3 & \cdots & A_h
\end{pmatrix},
\end{equation} where there is a diagonal identity matrix and 0 everywhere else, except for the last block row. 
Note that the effective information is all stored in the last block row of $A$ matrix. When we increase the number $h$ of delay coordinates, we increase length of this block row, thus increasing the linear dependency of the next state on the previous history state.

\subsection{Network Structure Detection}

Now we show that we can leverage the same $A$ matrix to detect the geometric network structure. In many cases it is desirable to detect network structure simply from traffic flow data. For example, if our data only contains real-time queue length data but lacks information on the geometric relationships of the queues, this structure detection technique would help us estimate the origin of the data and thus better predict future state or analysis stability.

In the case study we use the queue length data from the Montrose Rd \& E Jefferson St. intersection and the Montrose Pkwy \& E Jefferson St. intersection. We observe queue length data from 4 legs of each intersection and we would like to make inferences about the relative location of the two intersections. We stack the data of the two intersections as (NB1, SB1, WB1, EB1, NB2, SB2, WB2, EB2) so that the state vector $x_k \in \mathbb{R}^8$, where NB1 represents the NB leg of the first intersection we would like to detect. We choose $N=7200$ and $h=10$, meaning that we use $7200$ consecutive observations (2 hours) to form the data matrices $\tilde{X_1}$, $\tilde{X_2}$, as in \eqref{eqn:time-shifted-hankel}, after time-shifting them 10 times. We keep all the modes here as we use a rank truncation $\tilde{r} = 8 * 10$, the same as the dimension of $A$ matrix.

\begin{figure}[t]
	\centering
	\includegraphics[width=0.7\linewidth]{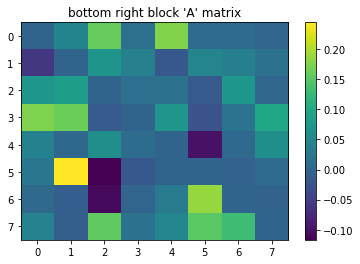}
	\caption{\label{fig:structure60} The bottom right block of A matrix structure from DMD with $10$ time-shifted data averaged over different dates of February 2018. The data that is binned at 60s. The entry with highest value on average in the matrix (the most yellow) represents the connection that SB1-SB2.}
\end{figure}


We run an experiment for each day of February 2018 using data from 0800-1000 and show the average of the results in Fig. \ref{fig:structure60}. We show the bottom right block A matrix from time-shifting DMD model, which is roughly equivalent to $A_h$ in \eqref{eq:arima} as shown in previous observation. Basically by investigating this matrix, we want to understand how $A_h x_t$ effects $x_{t+1}$. Note that the data is binned at 60s, which is about half the cycle time. This is important as there should be sufficient time for the flow out of one queue to enter the next connected queue. The entry with the highest value on average in Fig. \ref{fig:structure60} represents traffic flow from SB1 into SB2. This observation enables us to estimate that the first intersection is on the north of the second intersection. Thus we are able to detect that they are Montrose Rd \& E Jefferson St. intersection and Montrose Pkwy \& E Jefferson St. intersection respectively.

\subsection{Comparisons with Vector Autoregressions}
The last block row of \eqref{eqn:structure} can be written as
\begin{equation} \label{eq:arima}
x_{h+1} = \sum_{i=1}^{h} A_i x_i.
\end{equation}
Solving for a set of matrices $A_i$ such that \eqref{eq:arima} is approximately satisfied is exactly the problem considered in vector autoregression \cite{lutkepohl2005new}, which estimates relationships between the time series and their lagged values.


\begin{figure}[t]
	\centering
	\includegraphics[width=0.8\linewidth]{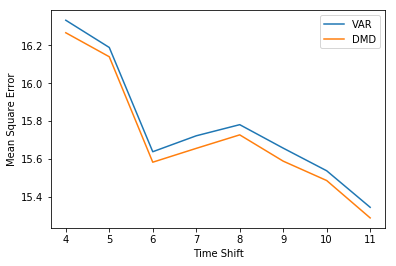}
	\caption{\label{fig:mse} Average Mean Square Error of prediction of queue length in the  next one minute using DMD and Vector Autoregressions (VAR).}
\end{figure}

To compare the prediction accuracy between DMD with time-shifted data and classic approaches to vector autoregression, we fit a vector auto regression model according to \eqref{eq:arima} by minimizing the least square error $|x_{h+1} - \sum_{i=1}^h A_i x_i|_2^2$ 
and compare to the result obtain using DMD as in \eqref{eqn:structure} for prediction. In both the DMD and vector autogression approach, we take $N=1200$ (20 minutes) observations and use the models to predict the next $60$ states (1 minute). We uniformly sample $20$ time windows in a day using queue length data at each second from the Montrose Rd \& E Jefferson St. intersection and compute the average of the mean square error with respect to different number $h$ of delay coordinates. As shown in Fig. \ref{fig:mse}, both methods tend to have better performance when increasing $h$. 

In general, as shown in Fig. \ref{fig:mse}, DMD outperforms vector autoregression for prediction when considering the same the number $h$ of delay coordinates. 


\section{Conclusion}
\label{sec:conclusion}
We have proposed using Koopman operator theory to study signalized traffic flow networks. The Koopman operator allows exactly representing a nonlinear dynamical system with a linear, but infinite dimensional, operator. This operator can be approximated with a finite-dimensional linear operator, \emph{i.e.}, matrix, and Dynamic Mode Decomposition provides a numerical algorithm for efficiently computing such an approximation using measured data. Moreover, this  approximation approach is model-free. We have demonstrated that this approach is applicable in a variety of applications for signalized traffic flow, such as inferring traffic control parameters, detecting and mitigating unstable queuing dynamics, and inferring structural relationships in a network of signalized intersections.

\section*{Acknowledgment}
The authors would like to thank Sensys Networks, Inc. for providing access to traffic flow data.


%





\ifCLASSOPTIONcaptionsoff
  \newpage
\fi





\bibliographystyle{IEEEtran}
\bibliography{IEEEabrv,ITS}

\vfill


\end{document}